\input amstex
\documentstyle{amsppt}
\nologo\magnification=1200
\overfullrule=0pt
\hsize = 5.5 true in\hoffset=.5in
\vsize = 9 true in

\loadeusm

\def\fM{{\frak M}}

\def\fmdtau{{\frak m_\tau}}

\def\frS{{\frak S}}
\def\frSds{{\frak S_*}}
\def\frSdsuz{{\frak S_*^0}}

\def\fX{{\frak X}}

\def\fW{{\frak W}}
\def\fWdz{{\fW_0}}

\def\sA{{\eusm A}}

\def\sB{{\eusm B}}

\def\sC{{\eusm C}}

\def\sE{{\eusm E}}

\def\sF{{\eusm F}}

\def\sL{{\eusm L}}

\def\sM{{\eusm M}}
\def\cM{{$\sM$ }}
\def\cMc{{$\sM$, }}
\def\cMp{{$\sM$. }}
\def\sMds{{\sM_*}}
\def\sMdp{{\sM_+}}
\def\tMdp{{\tM_+}}

\def\sMdsup{{\sM_*^+}}

\def\tM{{\widetilde {\sM}}}

\def\sN{{\eusm N}}

\def\sU{{\eusm U}}

\def\sR{{\eusm R}}

\def\C{{\Bbb C}}

\def\R{{\Bbb R}}

\def\Z{{\Bbb Z}}

\def\N{{\Bbb N}}

\def\a{{\alpha}}

\def\f{{\varphi}}
\def\p{{\psi}}

\def\Ad{{\text{\rm Ad}}}
\def\Adlr#1{{\text{\rm Ad}\lr{#1}}}

\def\edi{{e_i}}

\def\edn{{e_n}}

\def\edm{{e_m}}

\def\txti{{\text{\rm i}}}

\def\txm{{\text{\rm m}}}

\def\txmthth1{{\txm_{\thth1}}}

\def\thds{{\th_s}}

\def\th{{\theta}}

\def\tht{{\theta_T}}
\def\tht'{{\th_{T'}}}
\def\la{{\lambda}}
\def\ladz{{\la_0}}

\def\ladt{{\la_t}}

\def\rhods{{\rho_s}}
\def\rhodt{{\rho_t}}

\def\rhodg{{\rho_g}}

\def\ladz{{\la_0}}

\def\e{{\varepsilon}}

\def\Aut{{\text{\rm Aut}}}
\def\Autf'{{\Aut_\f'}}
\def\Autp'{{\Aut_\p'}}
\def\Cntp'{{\Cnt_\p'}}
\def\Cnt{{\text{\rm Cnt}}}

\def\Int{{\text{\rm Int}}}

\def\Intp'{{\Int_\p'}}

\def\mod{{\text{\ \rm mod\ }}}

\def\Autf{{\Aut_\f}}

\def\Autf'{{\Aut_\f'}}
\def\Proj{{\text{\rm Proj}}}

\def\two{{\rm I\!I}}
\def\twoone{{\rm I\!I$_1$}}
\def\threee{{\text{\rm I\!I\!I }}}

\def\threeone{{\rm I\!I\!I$_1$ }}
\def\threeonep{{\rm I\!I\!I$_1$. }}
\def\threeonec{{\rm I\!I\!I$_1$, }}

\def\threel{{\rm I\!I\!I$_{\lambda}$ }}

\def\three0{{\rm I\!I\!I$_0$}}

\def\bd{{\Bar\d}}

\def\g{{\gamma}}

\def\gdi{{\g_i}}

\def\gdone{{\g_1}}

\def\hdone{{h_1}}

\def\ggdi{{g_i}}

\def\ppdone{{p_1}}

\def\hdone{{h_1}}

\def\hdi{{h_i}}

\def\hdn{{h_n}}

\def\kdone{{k_1}}

\def\part{{\partial}}

\def\log{{\text{\rm log}}}

\def\sig{{\sigma}}

\def\sigf{{\sigma^{\f}}}
\def\sigft{{\sigma_t^{\f}}}
\def\sigfs{{\sigma_s^{\f}}}
\def\sigfr{{\sig^\f_r}}

\def\sigfditw{{\sig^\f_{\txti/2}}}

\def\sigfdmitw{{\sig^\f_{-\txti/2}}}

\def\sigfdmitwlr#1{{\sig^\f_{-\txti/2}\lr{#1}}}

\def\sigfdmi{{\sig^\f_{-\txti}}}

\def\sigufp{{\sig^{\f, \p}}}

\def\sigufpdt{{\sig^{\f, \p}_t}}

\def\sigufpdal{{\sig^{\f, \p}_\a}}

\def\wt{{semi-finite normal weight }}
\def\fwt{{faithful \wt}}
\def\botimes{{\Bar \otimes}}
\def\id{{\text{\rm id}}}

\def\r0{{\sR_0}}

\def\r01{{\sR_{0,1}}}

\def\botimes{{\overline \otimes}}
\def\wt{{\text\allowlinebreak{semi-finite normal weight}}}
\def\wts{{\text\allowlinebreak{semi-finite normal weights}}}
\def\fwt{{faithful \wt}}
\def\fwts{{faithful \wts }}

\def\vna{{\text\allowlinebreak{von Neumann algebra }}}

\def\vnap{{\text\allowlinebreak{von Neumann algebra. }}}
\def\vnasp{{\text\allowlinebreak{von Neumann algebras. }}}

\def\vnsa{{\text\allowlinebreak{von Neumann subalgebra}}}

\def\tB{{\text{\rm B}}}

\def\tH{{\text{\rm H}}}

\def\tZ{{\text{\rm Z}}}

\def\QED{{\hfill$\heartsuit$}}
\def\diamond{{
\centerline{-----
$\diamondsuit$$\diamondsuit$$\diamondsuit$
----- }\newpage}}

\def\inv{{^{-1}}}

\document

\def\etat{{{\eta_{T}}}}

\def\etat'{{\eta_{T'}}}

\def\_#1{{_{#1}}}
\def\^#1{{^{#1}}}

\def\(#1){{^{({#1})}}}
\def\ginv#1{{\g_{#1}^{-1}}}
\def\scirc{{\lower-.3ex\hbox{{$\scriptscriptstyle\circ$}}}}

\def\bracett'#1{{\left\{\!#1\right\}_{T'}}}
\def\bracet'm#1{{\bracett'{\txm\left({#1}\right)}}}

\def\brackett#1{{\left[{#1}\right]_{T}}}
\def\brackett'#1{{\left[{#1}\right]_{T'}}}

\def\rt'z{{\R/T'\Z}}

\def\tbth1{{\tB_\th^1}}

\def\tzth1{{\tZ_\th^1}}

\def\thth1{{\tH_\th^1}}

\def\pijj'{{\pi_{J, J'}}}
\def\pijj's{{\pi_{J, J'}^*}}

\def\brho{{\Bar\rho}}

\def\linfrt'z{{L^\infty(\R/T'\Z)}}
\def\linfr{{L^\infty(\R)}}
\def\linflr#1{{L^\infty\lr{#1}}}

\def\lone{{L^1}}
\def\lonelr#1{{L^1\lr{#1}}}
\def\loner{{L^1(\R)}}

\def\rt'z{{\R/T'\Z}}

\def\txd{{\text{\rm d}}}

\def\RN#1#2#3{{(\text{\rm D}#1:
\text{\rm D}#2)_{#3}}}
\def\dfn={{\overset{\text{\rm def}}\to=}}
\def\BBig(-{{\Big(\!\!-}}
\def\BBiglangle-{{\Big\langle\!\!-}}
\def\lrangle#1{{\left\langle{#1}\right\rangle}}
\def\:{{\text{\rm:}}}
\def\;{{\text{\rm;}}}
\def\M{{\text{\rm M}}}

\def\sumd#1{{\sum_{#1}}}

\def\Card{{\text{\rm Card}}}

\def\Linfr{{L^\infty(\R)}}

\def\intr{{\int_\R}}

\def\adi{{a_i}}

\def\bd#1{{\bd_{#1}}}

\def\bd#1{{b_{#1}}}

\def\kdone{{k_1}}
\def\kdn{{k_n}}

\def\kdn{{k_n}}

\def\r'dij{{r'_{i,j}}}
\def\r'dji{{r'_{j,i}}}

\def\pdone{{p_1}}

\def\qdone{{q_1}}

\def\qdi{{q_i}}

\def\pdi{{p_i}}

\def\edz{{e_0}}

\def\fdz{{f_0}}

\def\fdn{{f_n}}

\def\fdm{{f_m}}

\def\fdi{{f_i}}

\def\ffdi{{\f_i}}

\def\ffdone{{\f_1}}
\def\ffd2{{\f_2}}
\def\ppdone{{\p_1}}
\def\ppdtw{{\p_2}}

\def\ppdi{{\p_i}}

\def\udi{{u_i}}
\def\udius{{u_i^*}}

\def\udn{{u_n}}

\def\vdone{{v_1}}

\def\xdz{{x_0}}
\def\xdone{{x_1}}

\def\xdzus{{x_0^*}}

\def\ydone{{y_1}}

\def\bar xdone{{{\Bar x}_1}}

\def\xdi{{x_i}}

\def\xdius{{x_i^*}}

\def\ydone{{y_1}}
\def\xdn{{x_n}}

\def\ydn{{y_n}}

\def\ydi{{y_i}}

\def\bzdz\inv{{{\Bar z}_0^{-1}}}

\def\udone{{u_1}}

\def\udones{{u_1^*}}

\def\omegdxiz{{\omega_\xidz}}

\def\explr#1{{\exp\!\lr{#1}}}

\def\expums{{e^{-s}}}
\def\expus{{e^{s}}}

\def\lr#1{{\left({#1}\right)}}
\def\lrbrace#1{{\left\{#1\right\}}}\
\def\lrbracket#1{{\left[{#1}\right]}}

\def\sumdnonetoinf{{\sum_{n=1}^\infty}}

\def\sumd#1{{\sum_{#1}}}

\def\ot#1#2{{{#1}\otimes{#2}}}

\def\ffd#1{{f_{#1}}}

\def\ot{{\otimes}}

\def\lrbrace#1{{\left\{#1\right\}}}

\def\pidt{{\pi_T}}
\def\pidt'{{\pi_{T'}}}

\def\\dotpidk{{\dot\pi_{k}}}

\def\wdi{{w_i}}

\def\wdius{{w_i^*}}

\def\wdn{{w_n}}
\def\wdnus{{w_n^*}}

\def\sdt{{s_{T}}}
\def\sdt'{{s_{T'}}}

\def\ltwlr#1{{L^2\lr{#1}}}
\def\xidz{{\xi_0}}

\def\sMds{{\sM_*}}

\def\etadt'#1#2{{\eta_{T'}\lrbrace{#1, #2}}}
\def\ccyl{{cocycle\ }}

\def\aldt{{\a_t}}

\def\lrnorm#1{{\left\|{#1}\right\|}}

\def\lrabs#1{{\left|{#1}\right|}}

\def\lrnormsq#1{{\left\|{#1}\right\|^2}}

\def\expu#1{{e^{#1}}}

\def\rhods{{\rho_s}}
\def\rhodt{{\rho_t}}
\def\sAdf{{\sA_\f}}
\def\limd#1{{\lim_{#1}\ }}
\def\mudtau{{\mu_\tau}}
\def\sEdth{{\sE_\th}}
\def\Idth{{I_\th}}
\def\Idthlr#1{{I_\th}\lr{#1}}
\def\sAds{{\sA_*}}
\def\sAdsup{{\sA_*^+}}

\TagsOnRight 
\topmatter
\title{Noncommutative Integration }
\endtitle

\author{Masamichi Takesaki}
\address{Department of Mathematics, UCLA, P.O.Box 951555, Los Angeles, California 90095-1555. Mail Address: 3-10-39 Nankohdai, Izumi-ku, Sendai, 981-8003 Japan.}
\endaddress
\email{mt\@math.ucla.edu}
\endemail
\endauthor

\dedicatory{Dedicated to the Memory of Two Distinguished Operator Algebraists\: William B. Arveson and Gert K. Pedersen.}
\enddedicatory
\abstract{We will show that if $\sM$ is a factor, then for any pair $\f, \p\in\sMdsup$ of normal  positive linear functionals on $\sM$, the inequality: 
$$
\lrnorm{\f}\leq \lrnorm{\p}
$$ 
is equivalent to the fact that there exist a countable family $\lrbrace{\ffdi: i\in I}\subset \sMdsup$  in $\sMdsup$ and a family  $\lrbrace{\udi:  i\in I}\i\sM$ of partial isometries in \cM such that
$$
\f=\sumd{ i\in I} \ffdi,\quad \sumd{ i\in I} \udi{\ffdi}\udius\leq \p, \quad  \text{and} \quad \udius\udi=s\lr{\ffdi}, i\in I,
$$
where $s(\omega), \omega\in\sMdsup$, means the support projection of $\omega$. Furthermore, if $\lrnorm{\f}=\lrnorm{\p}$, then the equality replaces the inequality in the second statement. 
In the case that $\sM$ is not of type \threeonec the family of partial isometries can be replaced by a family of unitaries in \cMp 
One cannot expect to have this result in the usual integration thoery. To have a similar result, one needs to bring in some kind of non-commutativity. Let $\lrbrace{X, \mu}$ be a $\sig$-finite semifinite measure space and $G$ be an ergodic group of  automorphisms of $\linflr{X, \mu}$, then for a pair $f$ and $g$ of $\mu$-integrable positive functions on $X$, the inequality:
$$
\int_X f(x)\txd \mu(x)\leq \int_X g(x)\txd \mu(x)
$$
is equivalent to the existence of  a countable families $\lrbrace{\fdi:  i\in I}\subset L^1(X, \mu)$ of positive integrable functions and $\lrbrace{\gdi:  i\in I}$ in  $G$ such that
$$
f=\sumd{ i\in I} \fdi\quad\text{and}\quad \sumd{ i\in I} \gdi\lr{\fdi}\leq g,
$$
where the summation and inequality are all taken in the oredered Banach space $L^1(X, \mu)$ and the action of $G$ on $\lonelr{X, \mu}$ is defined through the duality between $\linflr{X, \mu}$ and $\lonelr{X, \mu}$, i.e.,
$$\aligned
 \lr{\g(f)}(x)&=f\lr{\g\inv x}\frac{\txd\mu\scirc \g\inv}{\txd\mu}(x), \quad f\in\lonelr{X, \mu}.
\endaligned
$$
 }
\endabstract
\leftheadtext{Noncommutative Integration}
\rightheadtext{Noncommutative Integration}

\endtopmatter
\document
\head{\S 1. Introduction}
\endhead
Regardless of commutativity, the integration of a positive element is the numerical value indicating the size of the quantity represented by the element. The one faces the following basic question:
$$\gathered
\text{What does two positive elements to record}\\
\text{ the same integration value mean? }
\endgathered
$$
Of course one cannot expect that two positive elements with the same integration value are isomorphic. In the classical integration theory, one cannot go further on this question. But in sharp contrast, in the non-commutative world, one can say that two positive elements with the same integration values are decomposed into the countable sum of two sequences of mutually isomorphic positive elements. This means that the non-commutative integration represents better the true meaning of integration than the classical commutative integration theory. Another important fact on this result is that the summation is taken over a countable set of objects. Otherwise, we are dealing with cardinality, which gives us very little room for analysis. This shows the distinguished position of the countability among infinities.  

\head{\S2. Preliminary, Noncommutative flow of weights}
\endhead
We will refer to either \cite{FT2} or \cite{ Tk2: Chapter XII  Section 6} for the basic facts on noncommutative flow of weights. But unfortunaely, \cite{Tk2: Exercise XII.6} contains a little inprecise statement, so we will present here the essence of that theory. 
We consider the translation flow $\lrbrace{\linfr, \R, \rho}$: 
$$
\lr{\rhodt f}(s)=f(s+t), \quad f\in\linfr, s, t\in \R.
$$
\proclaim{Lemma 2.1} If $\mu$ is a normal weight on $\sA=\linfr$ such that
$$\aligned
\mu\scirc \rhods(f)=\expums \mu(f), \quad f \in \sA_+;\\
0<\mu\lr{\fdz}<+\infty\quad \text{for some }\fdz\in\sA_+, 
\endaligned
$$
then the weight $\mu$ is a \fwt\ on $\sA$ such that
$$\gathered
C=\mu\lr{(-\infty, 0]}<+\infty;\\
\mu(f)=C\intr f(s)\expus \txd s, \quad f\in\linfr_+,
\endgathered
$$
where we view the normal weight $\mu$ as a measure on $\R$ absolutely continuous relative to the Lebesgue measure, but not necessarily semi-finite.
\endproclaim
\demo{Proof} Let $g$ be a continuous non-negative function with compact support on $\R$. Then we have 
$$\aligned
\mu\lr{\rhodg\lr{\fdz}}&=\mu\lr{\intr g(s)\rhods \lr{\fdz} \txd s}
=\intr g(s)\mu\lr{\rhods \lr{\fdz}} \txd s\\
&=\lr{\intr \expums g(s)\txd s }\mu\lr{\fdz},
\endaligned
$$
so that $0<\mu\lr{\rhodg\lr{\fdz}}<+\infty$ and 
$$\aligned
\lr{\intr g(s)\rhods \lr{\fdz}\txd s}(t)=\intr g(s)\fdz(s+t)\txd s=\lr{g*\fdz}(t).
\endaligned
$$
Hence $\rhodg\lr{\fdz}$ is continuous on $\R$ and takes a finite value on the normal weight $\mu$. Thus we may and do take a continuous positive function as $\fdz$ in the assumption of the lemma. So there are an interval $(a, b], a< b$, and a constant $C_1>0$ such that
$$
C_1\chi_{(a, b]}\leq \fdz \quad \text{and}\quad 0<\mu\lr{\chi_{(a, b]}}<+\infty.
$$
As $\rhods\lr{\chi_{(a, b]}}=\chi_{(a-s, b-s]}$, we have
$$
\mu\lr{(a-s, b-s]}=\expums \mu((a, b]) \quad \text{for every } s\in\R.
$$
From this, it follows that the measre $\mu$ takes a positive finite value on every finite interval and also that
$$\gathered\aligned
C=\mu\lr{(-\infty, 0]}&=\mu\lr{\bigcup_{n=1}^\infty(-n, -n+1]}\\
&=\sum_{n=1}^\infty \expu{-n}\mu((0,1])=\frac1{e-1}\mu((0, 1])<+\infty,\\
\mu\lr{(-\infty, s]}&=\mu\lr{\rho_{-s}\lr{ \chi_{(-\infty, 0]}}}=e^s\mu((-\infty, 0])<+\infty,\\
\endaligned\\
\txd \mu(s)=\mu\lr{(-\infty, 0]}e^s\txd s.
\endgathered
$$
This completes the proof. 
\QED
\enddemo
Fix a \vna $\sM$ and consider the associated noncommutative flow of weights $\lrbrace{\tM, \R, \tau, \th}$ to have
$$\gathered
\sM=\tM^\th, \quad \tau\scirc\thds=\expums \tau, \\
\sM'\cap \tM=\sC=\text{The Center of }\tM,\\
\lrbrace{\sC, \R, \th}=\text{The flow of weights on }\sM.
\endgathered
$$

Let $\fM$ be the algebra of all $\tau$-measurable densely defined closed operators affiliated to $\tM$. The following criteria for $\tau$-measurability is very useful and easy to manage: 
\roster
\item"{}" A densely defined closed positive operator $h$ affiliated to $\tM$ is {\bf$\pmb\tau$-measurable} if and only if there exists a positive number $\ladz>0$ such that
$$
\tau\lr{\chi_{[\ladz, +\infty)}\lr{h}}<+\infty.
$$
\endroster
We are going to write $E_\la=\chi_{[\la, +\infty)}\in \linfr$ for each $\la>0$. The algebra $\fM$ is graded by the noncommutative flow $\lrbrace{\thds: s\in\R}$ as seen below. 

Setting 
$$
\fM(\a)=\lrbrace{x\in\fM: \thds(x)=\expu{-\a s}x, \quad s\in\R}, \quad \a\in \C,
$$
we obtain the following: 
\roster
\item"i)" The original \vna \cM is the fixed point algebra of $\th$, which is exactly the equality:
$$
\sM=\fM(0).
$$
\item"ii)" For each $p>0$ we write 
$$
L^p(\sM)=\fM\lr{\frac 1p}.
$$
\item"iii)" The cases that $p=1$ and $p=2$ are of particular interest for us: 
$$\gathered
L^1(\sM)=\fM(1)=\lrbrace{x\in\fM\!: \thds(x)=\expums x, s\in\R},\\
\ltwlr{\sM}=\fM\lr{\frac12}=\lrbrace{x\in\fM: \thds(x)=e^{- s/2}x, s\in\R}.
\endgathered
$$
as it will be identified with the predual $\sMds$ and the standard form of \cMp
\item"iv)" If $\Re \a<0$, then 
$$
\fM(\a)=\lrbrace{0}.
$$ 
\endroster
We now consider the  operator valued weight $\Idth$ from $\tMdp$ to the extended positive cone $\widehat\sM_+$ of $\sM$: 
$$
\Idthlr{x}=\intr \thds(x)\txd s, \quad x\in \tMdp.
$$
As in the paper of Falcone - Takesaki, \cite{FT2}, we denote the element of $\fM$ corresponding to $\omega\in\sMds$ by $T(\omega)\in L^1(\sM)$ to avoid possible confusions and write
$$
\omega(1)=\int \txd T(\omega),
$$
which is defined to be the following value:
$$
\int \txd T(\omega)=\tau\lr{a^\frac12 T(\omega)a^\frac12}=\omega(1)
$$
for any $a\in\tMdp$ with $\Idth(a)=1$. The middle quantity $\tau\lr{a^\frac12 T(\omega)a^\frac12}$ does not depend on the choice of $a\in\tMdp$ with $\Idth(a)=1$ as shown in \cite{FT2: Theorem 3.12}.

\head{\S3. Comparison of Integrals}
\endhead

Let \cM be a fixed \vnap 
Fixing a pair $\f, \p\in \sMdsup$ with $p=s(\f), q=s(\p)\in \Proj(\sM)$, consider the one parameter group $\lrbrace{\sig^{\f, \p}_t: t\in\R}$ of isometries on $p\sM q$ defined by the following:
$$
\sig^{\f, \p}_t(x)=T(\f)^{\txti t}xT(\p)^{-\txti t}, \quad x\in \sM,
$$
which appears on the $(1, 2)$-corner of $\M_2(\C) \botimes \sM$ of the modular automorphism group $\sig^\rho$ of the balanced positive linear functional $\rho=\f\oplus \p$:
$$
\rho=\pmatrix \f&0\\0&\p\endpmatrix\in \lr{\M_2(\C) \botimes \sM}_*^+.
$$
We then consider the subspace $\sA\lr{\f, \p}$ of entire elements in $p\sM q$ relative to $\sigufp$, i.e., $\sA\lr{\f, \p}$ is the set of all those elements $x\in p\sM q$ such that 
the function: $t\in\R\mapsto \sigufpdt(x)\in \sM$ has entire extension to $\C$. We denote its value at $\a\in\C$ by $\sigufpdal(x)\in\sM$. Of particular interest to us is the value at the half imaginary unit:  $\pm\txti /2$ which is $\sig^{\f, \p}_{{\pm \txti}/2}(x)\in\sM$.

\proclaim{Lemma 3.1} If $x\in \sA(\f, \p), x\neq 0$, then the element $T\lr{\f}^\frac12 xT(\p)^\frac12 \in \lone\lr{\sM}$ has the property\:
$$\gathered
\lrabs{T(\f)^\frac12 xT(\p)^\frac12}\leq \lrnorm { \sig^{\f, \p}_{-\txti/2}(x)}\p,\\ 
\lrabs{\lr{T\lr{\f}^\frac12 xT(\p)^\frac12}^*}\leq \lrnorm { \sig^{\p, \f}_{-\txti/2}(x^*)}\f,\\
T(\f)^\frac12 xT(\p)^\frac12\neq 0.
\endgathered
$$
\endproclaim
\demo{Proof} We consider the path:
$$
t\in\R\mapsto \sigufpdt\lr{x}=T(\f)^{\txti t}xT(\p)^{-\txti t}\in \sA(\f, \p)\i p\sM q,
$$
which admits entire extension:
$$
\sig^{\f, \p}_z=T(\f)^{\txti z}xT(\p)^{-\txti z}\in \sA(\f, \p), \quad z\in\C.
$$
The evaluation at $-\txti /2$ gives
$$
\sig^{\f, \p}_{-\txti/2}(x)=T(\f)^{\frac12}xT(\p)^{-\frac12},
$$
so that we get
$$
T(\f)^\frac12 xT(\p)^\frac12=T(\f)^{\frac12}xT(\p)^{-\frac12}T(\p)=\sig^{\f, \p}_{-\txti/2}(x)T(\p)\in\lonelr{\sM}. 
$$
Thus we get the following easy conclusion:
$$\aligned
\lrabs{T(\f)^\frac12 xT(\p)^\frac12}&=\lrbracket{\lr{\sig^{\f, \p}_{-\txti/2}(x)T(\p)}^*\lr{\sig^{\f, \p}_{-\txti/2}(x)T(\p)}}^\frac12\\
&\leq \lrnorm{\sig^{\f, \p}_{-\txti/2}(x)}T(\p).
\endaligned
$$
The other inequality follows similarly. 

The non-triviality of the element $T(\f)^\frac12 xT(\p)^\frac12$ follows from the fact that $T(\p)$ is non-singular on the range of the projection $q$ and $T(\f)$ is also on the range of  $p$.
\QED
\enddemo

{\smc Definition 3.2.} A pair $\f, \p\in \sMdsup$ of normal positive linear functionals is said to be {\bf equivalent}  and written
$$
\f\sim \p
$$
if there exists a partial isometry $u\in \sM$ such that
$$
u^*u\geq s(\f),\quad uu^*\geq s(\p)\quad \text{and}\quad u\f u^*=\p,
$$
which automatically gives 
$$
\f=u^*\p u.
$$
If the above $u$ can be chosen to be unitary, then we say that $\f$ and $\p$ are {\bf unitarily conjugate} and write
$$
\f\equiv \p \quad \mod\ \Int(\sM).
$$

\proclaim{Lemma 3.3} If \cM be a factor, then 
every pair $\f, \p\in\sMdsup$ of non-zero normal positive linear functionals on \cM admits a pair $\ffdone, \ppdone\in\sMdsup$   such that
$$\gathered
0\neq \ffdone\leq \f, \quad 0\neq\ppdone\leq \p \quad \text{and}\quad \ffdone\sim \ppdone.
\endgathered
$$
In the case that if every non-zero normal positive linear functional $\omega\in\sMdsup, \omega\neq0$, majorizes a non-zero non-faithful positive linear fuctional $\omega_1\in\sMdsup, \omega_1\neq 0$, then the above $\ffdone$ and $\ppdone$ may be chosen to be unitarily conjugate, i.e.,
$$
\ffdone\equiv \ppdone \quad \mod\ {\Int(\sM)}.
$$ 
\endproclaim
\demo{Proof} Choose $x\in \sA(\f, \p), x\neq 0$ and set
$$
\rho=\frac1{\lrnorm{\sig^{\f, \p}_{-\txti/2}\lr{x}}}\f^\frac12
 x\p^\frac12\in \lonelr{\sM}=\sMds.
$$
Then with the polar decompostion:
$$
\rho=v\lrabs{\rho}
$$
we have
$$\gathered
0\neq\p_1=\lrabs{\rho}\leq \p\quad\text{and}\quad 0\neq \f_1=\lrabs{\rho^*}\leq \f,\\
v^*v=s\lr{ \p_1}, \quad vv^*=s\lr{ \f_1}\quad \text{and}\quad 
v\p_1v^*=\f_1, \quad v^*\f_1v=\p_1. 
\endgathered
$$
Setting $\udone=v^*$, we get the desired triplet $  \lrbrace{\ffdone,\ppdone, \udone}$  of the lemma. If the partial isometry $\udone$ admits a unitary extension $w$ in the sense that 
$$
w^*w=ww^*=1, \quad ws\lr{\ffdone}=\udone, 
$$
then the triplet $\lrbrace{\ffdone, \ppdone, w}$ is the required one in the latter claim. Thus if the projections $1-s\lr{\ffdone}$ and $1-s\lr{\ppdone}$ are equivalent in the projection lattice $\Proj(\sM)$, then the above $w$ exists and the last assertion on the unitary choice of $\udone$ follows. We split the proof according to the type of \cMp The case that \cM is finite has been taken care of by the above arguments. So we assume that \cM is infinite.
 
{\bf The case that \cM is semi-finite: } Let $\tau$ be a faithful semi-finite normal trace. Then $\ffdone$ and $\ppdone$ are of the following form:
$$\gathered
\ffdone(x)=\tau\lr{\hdone x} \quad\text{and}\quad \ppdone(x)=\tau\lr{\kdone x}, \quad x\in \sM,\\
\udone\hdone\udones=\kdone.
\endgathered
$$
Choose a spectral projection $e$ of $\hdone$ such that
$e\hdone\neq 0$ and $\tau(e)<+\infty$ and se $f=\udone e\udones$. Replacing the triplet $\lrbrace{\ffdone, \ppdone, \udone}$ by $\lrbrace{\ffdone e, \ppdone f, \udone e}$,  we can extend $\udone e$ to a unitary $w$, which makes the situation back to the already treated case. 

{\bf The case that \cM is purely infinite: }  Suppose \cM is purely infinite. In this case, every non-zero $\sig$-finite projections are equivalent and also the orthogonal complements of $\sig$-finite projections are equivalent in the case that \cM is not $\sig$-finite.  So if $s\lr{\ffdone}\neq 1$ and $s\lr{\ppdone}\neq 1$, then  we have
$$
1-s\lr{\ffdone}\sim 1-s\lr{\ppdone}.
$$
Thus we are back to the already treated case. Therefore the only remaining case is that eigther $s\lr{\ffdone}=1$ or $s\lr{\ppdone}\neq 1$ by symmetry. In this last case, the assumption on \cM guarantees the existence of a non-faithful $\omega\in\sMdsup, \omega\neq 0$ bounded by $\ffdone$, so that 
$e=s(\omega)\neq 1$. Replace $\ffdone$ by $\omega$ and set $\p_2=\udone\omega\udones$. Then we have $s\lr{\ppdtw}=\udone e\udones$ and 
$$
1-s(\ffdone)\sim 1-s\lr{\ppdtw}
$$
which allows us to extend $\udone s(\ffdone)$ to a unitary $w\in \sU(\sM)$ with $\udone s(\ffdone)=we$. This completes the proof of lemma. 
\QED
\enddemo

\proclaim{Lemma 3.4} Let $\lrbrace{\ffdi: i\in I}$ be a family of non-zero positive linear functionals on a \vna \cM such that there exists $\f\in\sMdsup$ which dominates all finite sums of $\f_i$, i.e., 
$$
\sumd{i\in J} \ffdi \leq \f\quad \text{for all finite subset }J\Subset I,
$$
then the family $\lrbrace{\ffdi: i\in I}$ is countable. 
\endproclaim
\demo{Proof} For each $n\in \N$, set
$$
I_n=\lrbrace{i\in I: \lrnorm{\ffdi}\geq \frac1n}.
$$
Then we have
$$
n\f(1)\geq \sumd{i\in I_n} n\ffdi(1)\geq {\Card\lr{I_n}},
$$
so that $\Card\lr{I_n}$ is finite. Since $I=\cup_{n\in\N}I_n$, we conclude that $I$ is countable.  

\QED
\enddemo
\proclaim{Theorem 3.6} \rm{\bf (Comparison of Positive Linear Functionals)}\ Let \cM be a factor.
 For a pair $\f, \p\in\sMdsup$ of non-zero normal positive linear functionals on \cMc the following statements are equivalent\: 
\roster
\item"i)" 
$$
\lrnorm{\f}=\f(1)\leq\p(1)=\lrnorm{\p}. 
$$
\item"ii)" There exist sequences $\lrbrace{\ffdi: i\in I}\i\sMdsup$, $\lrbrace{\ppdi: i\in I}\i\sMdsup$, $I\i \N$  such that
$$\gathered
\f=\sumd{i\in I} \ffdi, \qquad \sumd{i\in I} \ppdi\leq \p,\\
 \ffdi\sim\ppdi,\quad i\in I.
\endgathered
$$ 
\endroster
In the above equivalence, the equality of {\rm (i)} corresponds to that of {\rm (ii)}.
\endproclaim
\demo{Proof} Suppose that (i) holds. 
Let $\sF$ be the set of following three sequences:
$$\gathered
\Phi=\lrbrace{\ffdi: i\in I}\i \sMdsup, \quad \Psi=\lrbrace{\ppdi: i\in I}\i\sMdsup, \\
U=\lrbrace{\udi: i\in I}\i\sM
\endgathered
$$
such that
$$\gathered
\sumd{i\in } \ffdi\leq \f, \qquad \sumd{i\in I} \ppdi\leq \p,\\
 0\neq\udius\udi=s\lr{\ffdi}, \quad0\neq \udi\udius=s\lr{\ppdi}, \\
\udi\ffdi\udius=\ppdi \quad\udius\ppdi\udi=\ffdi, \quad i\in I.
\endgathered
$$
From Lemma 3.4 it follows that $\sF$ is an inductive set relative to the inclusion ordering. Hence it admits a maximal element
$\lrbrace{\Phi, \Psi, U}\in\sF$. The maximality and Lemma 2.3 implies that either
$$
\f=\sumd{i\in I} \ffdi\quad \text{or}\quad \p=\sumd{i\in I} \ppdi.
$$
If $\f\neq \sumd{i\in I} \ffdi$, then 
the equality $\p=\sumd{i\in I} \ppdi$ implies that
$$\aligned
\f(1)&>\sumd{i\in I} \ffdi(1)=\sumd{i\in I} \ppdi(1)=\p(1),
\endaligned
$$
which contradicts the assumption $\f(1)\leq \p(1)$. Hence we have
$$
\f=\sumd{i\in I} \ffdi\quad \text{and}\quad \p\geq \sumd{i\in I}\ppdi.
$$

Suppose that (ii) holds, i.e., there exists $\Phi=\lrbrace{\ffdi: i\in I}, \Psi=\lrbrace{\ppdi: i\in I}$ and $U=\lrbrace{\udi: i\in I}$ which satisfies the requirements in (ii). Since $\udi$ is an isometry from the range of $\pdi=s\lr{\ffdi}$ to that of $\qdi=s\lr{\ppdi}$, we have $\lrnorm{\ffdi}=\lrnorm{\ppdi}, {i\in I}$. 
Then we get
$$\aligned
\f(1)&=\sumd{i\in I} \ffdi(1) =\sumd{i\in I} \lrnorm{\ffdi}=
\sumd{i\in I} \lrnorm{\ppdi}=\sumd{i\in I}\ppdi(1)\\
&\leq \p(1).
\endaligned
$$
This completes the proof.
\QED
\enddemo
\vskip.2cm
\noindent
{\smc Definition 3.7.} A positive linear functional $\f$ on a \vna \cM is said to be {\bf super} faithful if every non-zero positive linear functional $\p$ dominated by $\f$ is faithful. 

\vskip.2cm
\noindent
{\smc Remark 3.8.} If $\f$ is a super faithful state on a \vna $\sM$, then $\f$ is automatically a normal faithful positive linear functional and $\sM_\f=\C$, consequently $\sM$ is a factor of type \threeonep

\proclaim{Corollary 3.9} {\rm i)}If \cM is a factor which does not admits a super faithful state, then the equivalence $\ffdi\sim \ppdi$ in the condition {\rm(ii)} can be replaced by the unitary conjugacy\: $\ffdi\equiv \ppdi \mod \Int(\sM), {i\in I}$. 

{\rm ii)} If the pair $\f, \p\in\sMdsup$ are both non-faithful instead, then the equivalence $\ffdi\sim\ppdi$ can be replaced by 
$\ffdi \equiv \ppdi \mod \Int(\sM)$.

{\rm iii)} If the pair $\f, \p\in\sMdsup$ are both  super faithful, then the equivalence $\ffdi\sim\ppdi$ is replaced by 
$\ffdi \equiv \ppdi \mod \Int(\sM)$.
\endproclaim
\demo{Proof} This follows from the fact that the equivlance of non-faithful normal positive linear functionals can be implemented by a unitary element in \cMp
\enddemo
 \head{\S4. Commutative Case}
\endhead
Let $\sA$ be an abelian \vnap 
In this case, as it stands, one cannot compare a pair of normal positive linear functionals beyond the absolue continuity. We need  a device to move around elements of $\sA$. So let $G$ be a group of automorphisms of $\sA$, i.e., $G$ is a subgroup of $\Aut(\sA)$. For each member $\g\in G$, we consder the action of $\g$ on the predual $\sAds$ as follows:
$$
\lrangle{x, \g(\f)}=\lrangle{\ginv{} (x), \f}, \quad x\in\sA, \f\in\sAds.
$$
We write 
$$
\f\equiv \p\quad \mod G
$$
if there exists $\g\in G$ such that $\p=\g(\f)$.
\proclaim{Proposition 4.1} If $\sA$ is an abelian \vna equipped with an ergodic group $G$ of automorphisms, then for every pair $\f, \p\in \sAdsup$ of normal positive linear functionals the following two staatements are equivalent\:
\roster
\item"i)"
$$
\lrnorm{\f}\leq \lrnorm {\p}.
$$
\item"ii)" There exist families $\lrbrace{\ffdi: i\in I}$ and $\lrbrace{\ppdi: i\in I}$ of normal positive linear functionals on $\sA$ such that
$$\gathered
\f=\sumd{i\in I} \ffdi, \quad \sumd{i\in I}\ppdi \leq \p;\\
\ffdi \equiv  \ppdi\quad \mod G, \quad i\in I.
\endgathered
$$
\endroster
\endproclaim
\demo{Proof} First we remark that the commutativity of $\sA$ entails the lattice property of both the self-adjoint part of   $\sA$ and of the self-adjoint part of  its predual $\sAds$.
From the discussion in the last section, to prove the theorem it is enough to show that for every pair $\f, \p\in\sAdsup$ of non-zero normal positive linear functionals there exists a pair 
$\ffdone, \ppdone\in \sAdsup$ such that 
$$\gathered
0\neq \ffdone \leq \f, \quad 0\neq \ppdone \leq \p, \\
\ffdone\equiv \ppdone \quad \mod G.
\endgathered
$$
To this end, set
$$
p=s(\f), \quad q=s(\p).
$$
Since $\f\neq 0$ and $\p\neq 0$, we have $p\neq 0$ and $q\neq 0$ as well. Hence the ergodicity of $G$ implies the existence of $\gdone \in G$  such that 
$$
\gdone(p)q\neq 0.
$$
This means that $\gdone(\f)q\neq 0$, so that 
$\gdone(\f)\wedge \p=\ppdone \neq 0$. Setting 
$$
\ffdone=\ginv{1}\lr{\ppdone},
$$
we obtaind a pair $\ffdone, \ppdone\in\sAdsup$ with 
$$
0\neq \ffdone\leq \f, \quad 0\neq \ppdone\leq \p, \quad \ffdone \equiv \ppdone \mod G.
$$
This completes the proof. \QED
\enddemo

Application of the proposition yields the following fact which can be stated more general form such as the integration over a locally compact group. We just state here a special case which should be taught in the class on the Lebesgue integration.

\proclaim{Corollary 4.2} Let $\lonelr{\R^n}$ be the Banach space of all integrable functions on the vector space $\R^n$ relative to the Lebesgue measure. For a pair $f, g\in \lonelr{\R^n}_+$ of positive integrable functions, the following two conditions are equivalent\:
\roster
\item"i)" 
$$
\int_{\R^n}f(x)\txd x\leq \int_{\R^n } g(x)\txd x.
$$
\item"ii)" There exist countable families, $\lrbrace{\fdi: i\in I}, \lrbrace{\ggdi: i\in I}\i \lonelr{\R^n}_+$ and $\lrbrace{\adi: i\in I}\i \R^n$ such that
$$\gathered
f=\sumd{i\in I} \fdi, \quad \sumd{i\in I}\ggdi \leq g \\ \text{and}\\ 
\ggdi(x)=\fdi\lr{x+\adi}\ \text{for almost every } x\in\R^n.
\endgathered
$$
Here the summation is taken relative to  the convergence in the Banach space $\lonelr{\R^n}$.
\endroster
Here the equality of {\rm(i)} corresponds to that of {\rm (ii)}.
\endproclaim

\head{\S5. Commutativity of Normal Positive Linear Functionals}
\endhead
Fix a factor  \cM and a pair $\f, \p\in \sMdsup$ with $\lrnorm{\f}\leq \lrnorm{\p}$. Then we have decompostion:
$$
\f=\sumd{i\in I}\ffdi, \quad \sumd{i\in I}\ppdi \leq \p, \quad 
\ffdi\sim \ppdi,\quad  i\in I.
$$
We are going to discuss the commutativity of the families $\lrbrace{\ffdi: i\in I}$ and $\lrbrace{\ppdi, i\in I}$. To this end, we remind ourselves the following fact:
the commutativity of a pair $\f, \p\in\sMdsup$ of normal positive linear functionals was first introduced in \cite{Tk1} in the following form: 
\vskip.2cm
\noindent
{\smc Definition 5.1.} A pair $\omega_1, \omega_2\in\sMdsup$ is said to {\bf commute} if 
$$
\lrabs{\omega_1+\txti \omega_2}=\lrabs{\omega_1-\txti \omega_2}.
$$

In the case that both functionals are faithful, it is shown in \cite{Tk1} that their commutativity  is equivalent to the invariance of one relative to the modular automorphism group of the other. For the general pair $\f, \p\in\sMdsup$, we don't have any tool to attack the commutativity question. So we restrict ourselves to the special case that $\f$ and $\p$ are both factoring through a maximal abelian subalgebra $\sA$ of \cM in the sense that 
$$
\f=\f\scirc \sE_\sA \quad\text{and}\quad \p=\p\scirc \sE_\sA,
$$
where $\sE_\sA$ means the $\sA$-valued normal conditional expection. In general, $\sE_\sA$ does not exist. For example there is no normal conditional expectation from $\sL\lr{L^2(\R)}$ to $\linfr$. But if it does exist, then it is unique. 

\proclaim{Proposition 5.2} Let \cM be a factor and $\sA$ be a maximal abelian subalgebra of \cMp 
If $\sA$ is semi-regular and the range of normal conditional expectation $\sE_\sA$, then for a pair $\f, \p \in \sMdsup$ such that 
$$
\f=\f\scirc \sE_\sA, \quad \p= \p\scirc \sE_\sA,
$$
the inequality 
$$
\lrnorm{\f}\leq \lrnorm{\p}
$$
is equivalent to the existence of the decompostion\: 
$$\gathered
\f=\sumd{i\in I}\ffdi, \quad \sumd{i\in I}\ppdi\leq \p,\\
\ffdi \equiv \ppdi \mod \sN(\sA),\\
\ffdi=\ffdi\scirc \sE_\sA, \quad \ppdi=\ppdi\scirc \sE_\sA, \quad i\in I,
\endgathered
$$
where $\sN(\sA)=\lrbrace{u\in\sU(\sM): u\sA u^*=\sA}$ is the normalizer of $\sA$ in \cMp
 \endproclaim
Before the proof, we observe that the invariance $\f=\f\scirc \sE_\sA$ is equivalent to the inclusion:
$$
\sA\i \sM_\f.
$$
\demo{Proof} First we observe 
$$
u\sE_\sA(x)u^*=\sE_\sA\lr{uxu^*}, \quad u\in \sN(\sA).
$$
Then with $G=\lrbrace{\Ad(u): u\in \sN(\sA)}$, $G$ acts on $\sA$ ergodically by the semi-regularity assumption on $\sA$. Hence Proposion 4.1 implies that there exist families $\lrbrace{\Bar\ffdi: i\in I}\i \sAdsup$ and $\lrbrace{\Bar\ppdi: i\in I}\i \sAdsup$ such that
$$
\f|_\sA=\sumd{i\in I} \Bar\ffdi, \quad \sumd{i\in I}\Bar\ppdi \leq \p|_\sA, \quad \Bar\ffdi\equiv \Bar\ppdi\quad \mod G.
$$ 
Setting $\ffdi=\Bar\ffdi\scirc \sE_\sA$ and $\ppdi=\Bar\ppdi\scirc \sE_\sA$, we get for every $x\in\sM_+$, 
$$\aligned
\sumd{i\in I} \ffdi(x)&=\sumd{i\in I} \Bar\ffdi\lr{ \sE_\sA(x)}
=\lr{\sumd{i\in I}\Bar\ffdi}\lr{\sE_\sA(x)}\\
&=\f\lr{\sE_\sA(x)}=\f(x);\\
\sumd{i\in I} \ppdi(x)&=\sumd{i\in I} \Bar\ppdi\lr{ \sE_\sA(x)}
=\lr{\sumd{i\in I}\Bar\ppdi}\lr{\sE_\sA(x)}\\
&\leq\p\lr{\sE_\sA(x)}=\p(x).
\endaligned
$$
If $\udi\in\sN(\sA)$ gives $\Adlr{\udi}|_\sA\lr{\Bar\ffdi}=\Bar\ppdi, i\in I$, then we have for each $x\in\sM$
$$\aligned
\ffdi\lr{\udius x\udi}&=\Bar\ffdi\lr{\sE_\sA\lr{\udius x \udi}}
=\Bar\ffdi\lr{\udius \sE_\sA(x)\udi}=\Bar\ffdi\scirc \Adlr{\udi}^{-1}\lr{\sE_\sA(x)} \\
&=\Bar\ppdi\lr{\sE_\sA(x)}=\ppdi(x).
\endaligned
$$
Consequently we get
$$
\ffdi\equiv \ppdi\quad \mod G, \quad i\in I.
$$
This completes the proof.
\QED
\enddemo

\head{\S6. Concluding Remark}
\endhead
Throughout the paper, we only consider the factor case. The generalization to the non-factor case is very much the same as the comparison theory of projections in the general frame work of \vnasp However the point of this paper is that the non-commutative theory of integration gives the natural answer about the question concerning the meaning of the same values on integrations. In the case of semi-finite factors, the work of Kadison and Pedersen, \cite{KP}, on the additivity of a trace gives the same result. However  the motivations of their work and this work are quite different. They are very much concerned about the natural proof of the additivity property of the trace which comes from the comparison of projections. In other words, their theory can be viewed as the one about the measure theory, whilst our work is more concerned with the result of integration. Technically, their work is more demanding as they don't assume the existence of a semi-finite normal trace on the base \vnap
Indeed, the result in the case of factors of type \threee is unexpected. Also the seek of natural answer brought about the new question on the existence of a super faithful state which was never considered before. The author has been unable to exclude the existence of a super faithful state so far. The author would like to leave the existence question of a super faithful state as a challenge for operator algebraists.

\enddocument 

\diamond
Is it possible to have $\f\in\sMdsup$ such that every non-zero
$\p\leq \f$ is faithful? If this occurs, then $\sM_\f=\C$. It means that $\sigfdmitw(a)$ for every non-zero $a\in\sA(\f)$ is non-singular. Thus
$$
\f^\frac12 h\f^{-\frac12} \text{ is non-singular}, \quad \forall h\in\sMdp, h\neq 0.
$$
$$
\lr{e^{s/2}\sigfditw+e^{-s/2}\sigfdmitw}^{-1} (h)=\intr \frac{e^{-\txti st}}{e^{\pi s} +e^{-\pi t}}\sigft(h)\txd t
$$
$$\aligned
a&=\lr{e^{s/2}\sigfditw+e^{-s/2}\sigfdmitw}\lr{\intr \frac{e^{-\txti st}}{e^{\pi s} +e^{-\pi t}}\sigft(h)\txd t}
\endaligned
$$
Suppose 
$$\gathered
h=\sig^\f_f(k)=\intr f(r)\sigfr(k)\txd r\\
\sigft(h)=\intr f(r-t)\sigfr(k)\txd r=\sig^\f_{\ladt f}(k).
\endgathered
$$
$$\aligned
h&=\lr{e^{s/2}\sigfditw+e^{-s/2}\sigfdmitw}\lr{\intr \frac{e^{-\txti st}}{e^{\pi s} +e^{-\pi t}}\sig^\f_{\ladt f}\txd t}(k)\\
&=\lr{e^{s/2}\sigfditw+e^{-s/2}\sigfdmitw}\lr{\sig^\f_{\lr{e^{s/2}\la_{\txti/2}+e^{-s/2}\la_{-\txti/2}}^{-1}f}}^{-1}(k).
\endaligned
$$

$$\aligned
\intr \frac {e^{-\txti st}}{e^{\pi t}+e^{-\pi t}}\sigft(k)\txd t
&=\lr{e^{s/2}\sigfditw+e^{-s/2}\sigfdmitw}^{-1}(k)\\
&=e^{-s/2}\sigfdmitw\lr{\id + e^{-s}\sigfdmi}^{-1}(k)
\endaligned
$$
$$\aligned
\sigfdmitw(h)^*&\sigfdmitw(h)=\sigfditw\lr{h}\sigfdmitw(h)
\ \text{ is non-singular for every } h\in\sMdp.
\endaligned
$$
Furthermore, 
$$
T(\f)^\frac12hT(\f)^\frac12 \text{ is non-singular for every }h\in\sMdp.
$$
The set
$$
K=\lrbrace{\omega\in\sMdsup: 0\leq \omega \leq \f }\i\sMdsup
$$
is a weakly compact convex subset of $\sMdsup$, which is isomorphic to the $\sig$-weakly compact set:
$$
S_+=\lrbrace{h\in\sM: 0\leq h\leq 1},
$$
under the map:
$$
a\in S_+\mapsto \f^\frac12 a\f^\frac12\in K.
$$
The extreme boundary of  $S_+$  is exactly the set of projections:
$$
\Proj(\sM)=\part_e \lr{S_+}.
$$
Hence 
$$
\part_e(K)= \lrbrace{\f^\frac12 p\f^\frac12: p\in\Proj(\sM)}.
$$

$$
T(\f)^\frac12\xi=pT(\f)^\frac12\xi+(1-p)T(\f)^\frac12\xi
$$

$$\gathered
0\neq pT(\f)^\frac12\xi \neq T(\f)^\frac12\xi, \quad\forall p\in\Proj(\sM), \xi\neq 0, p\neq 0.
\endgathered
$$
Set
$$
\omega=\f^\frac12p\f^\frac12\leq \f.
$$
Then we have
$$\gathered
\omega^\frac12=\lr{\f^\frac12p\f^\frac12}^\frac12
=\lr{\sigfdmitwlr{p}\f}^\frac12=a\f^\frac12=\f^\frac12a^*\\
\omega=\f^\frac12a^*a\f^\frac12=\f^\frac12p\f^\frac12\\
p=a^*a, \quad aa^*=1.
\endgathered
$$
Suppose 

Suppose that the modular automorphism group $\sigf$ of a faithful normal state $\f\in\sMdsup$ is asymptotically abelian in the sense that 
$$
\lrnorm{\sigft(a)\omega-\omega\sigft(a)}\longrightarrow 0 \ \text{as } t\to \pm\infty.
$$

\proclaim{Lemma 2.2} If \cM is a factor, then for every non-zero element $\f\in\sMdsup$, there exists a pair of non-zero positive element $a, b\in \Proj(\sM)$  such that 
$$\gathered
aT\lr{\f}^\frac12 b=0\quad \text{but}\quad T(\f)^\frac12 b\neq 0.
\endgathered
$$
\endproclaim
\demo{Proof} With $p=s(\f)\in\Proj(\sM)$ the support projection of $\f$, we restrit our attention to the reduced factor $\sM_p$ so that we may and do assume that $\f$ is faithful. Suppose that our assertion fails, i.e., for evey pair of non-zero positive elements $a, b\in \sMdp$ we have
$$
a\f^\frac12 b\neq 0.
$$
This means that $\f^\frac12 h\f^\frac12$ is non-singular for every $h\in\sMdp$. If $h\in \sMdp\cap \sA(\f)$, then we have
$$
\f^\frac12 h\f^\frac12=\sigfdmitw(h)\f,
$$
and consequently the non-singularity of $\sigfdmitw(h)\f$.

This means that the modular automorphism group $\sigf$ has the property that 
$$
\bigvee_{t\in\R} \sigft(p)=1 \quad \text{for every non-zero } p\in \Proj(\sM).
$$
Consequently, the modular automorphism group $\sigf$ is ergodic, i.e., $\sM_\f=\C$. Furthermore, every non-zero positive linear functional bounded by $\f$ is faithful.

Since $\thds(a)=\expums a$, the range projection $p$ of $a$ is left invariant under $\thds, s\in\R$, so that it is a projection of $\sM$. Let $a=vk$ be the polar decompostion of $a$ with $k=\lrabs{a}$ and set $h=vkv^*=\lrabs{a^*}$. The phase operator $v$ belongs to \cM and $v^*v=q\in\Proj(\sM), vv^*=p$.
Let $\f$ and $\p$ be the elements in $\sMdsup$ such that 
$$
T(\f)=h\quad \text{and}\quad T(\p)=k.
$$
Then we have $s(\f)=p$ and $s(\p)=q$ and 
$$
\f=\p\scirc \Ad(v^*)\quad \text{on}\quad \sM_p\quad \text{and}\quad \p=\f\scirc \Ad(v)\quad\text{on}\quad \sM_q.
$$
\enddemo

$$\aligned
\omegdxiz(x)&=\lr{x\xidz\mid\xidz}\\
0=x\f^\frac12 a&=\f^\frac12\sig^\f_{\txti /2}\lr{x}a\\
0&=x\sig^\f_{-\txti/2}(a).
\endaligned
$$
Find  $a$ in \cM such that
$$\aligned
\ltwlr{\sM}&\neq\lrbracket{\f^\frac12 a\sM} =p\ltwlr{\sM}
\endaligned
$$

\proclaim{Lemma 2.3} If \cM is a factor, then for any  pair $\f, \p\in\sMdsup$ of non-zero elements there exist a pair $\f_1, \p_1$ of non-zero elments in $\sMds$ and a unitary $\udone\in\sU(\sM)$ such that  $0\neq\f_1\leq \f$, $0\neq \p_1\leq \p$  
$$
\f_1=\p_1\scirc \Adlr{\udone}.
$$
\endproclaim
\demo{Proof} In the last lemma, set
$$
a=\frac1{\lrnorm { \sig^{\f, \p}_{-\txti/2}(x)}}T(\f)^\frac12 xT(\p)^\frac12\in \lonelr{\sM}.
$$
Then by the last lemma, we have
$$
k=\lrabs{a}\leq \p\quad \text{and}\quad h=\lrabs{a^*}\leq \f.
$$
Let $a=vk$ be the polar decomposition of $a\in\lonelr{\sM}$, so that $v^*v=q$ and $vv^*=p$. Then we have $v^*hv=k$ and $vkv^*=h$. We split the case according to the equivalence of the orthogonal complement of  subprojections $\pdone$ and $\qdone$ of $p$ and $q$ respectively which commute with $h$ and $k$ respectively. 

{\bf The Case that $\pmb{1-\pdone\sim1-\qdone:}$\quad }
The case that  $1-p\sim 1-q$ is easy. As  there exists a partial isometry $w\in\sM$ such that 
$$
w^*w=1-q \quad\text{and}\quad ww^*=1-p.
$$
Setting
$$
\udone=v+w,
$$
we get a unitary $\udone\in\sU(\sM)$ such that 
$$
\udones k\udone=h\quad \text{and}\quad \udone h\udones =k.
$$
With $\f_1, \p_1\in\sMdsup$ such that $T\lr{\f_1}=h$ and $T\lr{\p_1}=k$,  we get the desired conclusion:
$$
\p_1=\f_1\scirc \Ad\lr{\udone}, \quad 0\neq \f_1\leq \f, \quad 0\neq \p_1\leq \p.
$$

{\bf The Case that $\pmb{1-p\not\sim1-q:}$\quad } In this case, the projections $p$ and $q$ are infinite. If $k$ admits a  finite spectal projection $\qdone$, then 
$\vdone=v\qdone$ admits a unitary extension $\udone$.  Setting $\kdone=k\qdone\neq0$, $\hdone=v\kdone v^*$ and choosing $\f_1$ and $\p_1$ as above,  we can conclude the proof. Thus the only case that remains is the case that every non-zero spectral projecion of $k$ is infinite.

\enddemo

\diamond

Fix a \vna $\sM$. Associated with $\sM$ is the noncommutative flow of weights $\lrbrace{\sM, \R, \th}$ so that 
$$\gathered
\sM=\tM^\th, \\
 \tau\scirc \thds=\expums \tau,\\
\sM'\cap\tM=\sC=\text{The center of }\tM.
\endgathered
$$
The \vna $\tM$ is called the core of $\sM$ and generated by
$\sM$ together wiht one parameter unitary groups $\lrbrace{\f^{\txti t}: t\in \R, \f\in \fW_0}$, where $\fW_0$ is the set of all \fwts on \cMp When we view $\f, \f\in \fW$, as an operator affiliated with $\tM$, we write it $h(\f)$ to avoid  possible confusion. We know that $\f\in \fWdz$ is finite if and only if $h(\f)$ is $\tau$-measurable in the sense that
$$
\limd{\la\to+\infty} \tau\lr{E_\la(h(\f))}=0
$$
which is equivalent to the fact:
$$
\tau\lr{\lr{E_\la(h(\f))}}<+\infty 
$$
for some large $\la>0$. It follows from the work of Falcone and Takesaki, \cite{FT2}, that the predual $\sMds$ is naturally identified with 
$$\aligned
\lone\lr{\sM}&=\lrbrace{x\in\fM: \thds(x)=\expums x, s\in \R}\\
&=\lrangle{xh(\f): x\in\sM, \f\in\sMdsup}
\endaligned
$$
where $\fM$ is the algebra of densely defined $\tau$-measurable closed operators and the notation $\lrangle{S}$ means the linear span of $S$. 

If $x\in\lone\lr{\sM}$, then the polar decomposition:
$$
x=u\lrabs{x}
$$
gives that $\lrabs{x}\in L^1(\sM)$. Suppose $h\in \lone(\sM)_+$, i.e.,
$$
\limd{\la\to+\infty}\tau(E_\la(h))=0 \quad \text{and}\quad
\thds(h)=\expums h, \quad s\in\R.
$$
With $p=s(h)\in\tM^\th=\sM$, restricting our attention to the reduced \vna $\sM_p$, we may and do assume that p=1.
$$
\thds\lr{h^{\txti t}}=\expu{-\txti st}h^{\txti t}, \quad s, t\in\R.
$$
With $\f\in\fWdz$, we have
$$
u_t=h^{\txti t}h(\f)^{-\txti t}\in \sU\lr{\tM^\th}=\sU(\sM), \quad t\in \R.
$$
Then the one parameter family $\lrbrace{u_t: t\in\R}$ is a \ccyl relative to $\sigf$:
$$\gathered
u_{s+t}=u_s\sigfs\lr{u_t}.
\endgathered
$$
Hence there exists $\p\in\fWdz$ such that
$$\gathered
u_t=\RN{\p}{\f}{t}=h(\p)^{\txti t}h(\f)^{-\txti t},\\
 h^{\txti t}=h(\p)^{\txti t},\\
\therefore\quad h=h(\p),
\endgathered\quad t\in \R. 
$$

Let $\sA$ be the abelian \vnsa\ generated by $h$.
Then the covariant system $\lrbrace{\sA, \R, \th}$ is identified with $\loner$ equipped with the translation $\rho$, i.e., the element $h$ is represented by the function:
$$
h(s)=\expums, \quad s\in \R,
$$
and with $\rho$ the translation flow:
$$
(\rhods f)(t)=f(t+s), \quad f\in\linfr, s, t\in\R,
$$
we have
$$
\lrbrace{\sA, \R, \th}\cong \lrbrace{\linfr, \R, \rho}.
$$
So, we will identify $\lrbrace{\sA, \R, \th}$ with $\lrbrace{\linfr, \R, \rho}$.  
The $\tau$-measurability of $h$ means that
$$
\tau\lr{E_\ladz(h)}<+\infty \quad \text{for some }\ladz>0.
$$
Since 
$$
E_\la(h)=\chi_{(-\infty, -\log \la]},
$$
the measure $\mudtau$ corresponding to $\tau$ gives a finite value to the half line $(-\infty, -\log\ladz]$:
$$
\mudtau\lr{(-\infty, -\log\ladz]}<+\infty.
$$
Furthermore, the translation:
$$\aligned
(\rhods E_\la)(t)&=(\rhods\chi_{(-\infty, \log\la]})(t)
=\chi_{(-\infty, \log\la]}(t+s)\\
&=\chi_{(-\infty, \log\la-s]}(t)
\endaligned
$$
gives
$$\aligned
\mudtau\lr{(-\infty, \log\ladz-s]}&=\tau\lr{\thds\lr{E_\ladz(h)}}
=\expums \tau\lr{E_\ladz(h)}\\
&=\expums \mudtau\lr{(-\infty, \log\ladz]}<+\infty;\\
\mudtau\lr{(-\infty, \log\ladz+s]}&=\expus \mudtau\lr{(-\infty, \log\ladz]}
\endaligned
$$
Consequently, we get
$$
\mudtau\lr{(0, s]}=\lr{\expus-1}\mudtau\lr{(-\infty, 0]}.
$$
Thus we get
$$
\txd \mudtau(s)={\mudtau\lr{(-\infty, 0]}}\expus \txd s.
$$
Recall 
$$
\tau h(\p)=\widehat\p=\p\scirc \sEdth, 
$$
so that the \fwt\ $\widehat\p$ on $\sA=\linfr$ is precisely  the 
integration against the measure: $\mudtau((-\infty, 0])\times$the Lebesgue measure.

We are now going to examine the finiteness of the \wt\ $\p$.
Define the operator valued \wt\ $\sEdth$ by the integration:
$$
\sEdth(a)=\intr \thds(a)\txd s, \quad a\in \tM_+.
$$
Since $\th$ leaves $\sA$ globally invariant, $\sEdth(\sA)\i \sA\cap \sM=\C$,  For each $f\in\linfr\cap \loner_+$ we have
$$\gathered
\sEdth(f)=\intr \thds\lr{f}\txd s
=\intr {f(s)}\txd s\in \R_+,\\
\tau\lr{h(\p)f}={\mudtau\lr{(-\infty, 0]}}\intr f(s)\txd s <+\infty,\\
\lrnorm{\p}=\p(1)=\widehat\p(f) = \mudtau\lr{(-\infty, 0]}<+\infty
\endgathered
$$
with $f\in \loner_+\cap \linfr$ such that
$$
\intr f(s)\txd s =1.
$$

We write
$$
\int x\f=\lrangle{1, x\f}=\lrangle{x, \f}=\f(x), \quad x\in\sM, \f\in\sMdsup.
$$

\proclaim{Main Theorem} Let $\sM$ be a factor. For a pair $\f, \p\in \sMdsup$ of normal positive linear functionals on $\sM$, the following conditions are equivalent\:
\roster
\item "i)" 
$$
\int \f=\f(1)\leq \int \p=\p(1),
$$
\item"ii)" There exist pair $\lrbrace{\f_n: n\in\N},\lrbrace{\p_n:\in\N}$ of  sequences of normal positive linear functions on $\sM$  such that
$$
\f=\sumdnonetoinf  \f_n\quad\text{and}\quad \quad \sumdnonetoinf \p_n\leq \p, \quad \f_n\equiv \p_n\quad \mod\lr{\Int(\sM)},
$$
where the notation $\f\equiv \p\mod (\Int(\sM))$ means the existence of a unitary $u\in\sU(\sM)$ such that $\p=\f\scirc \Ad(u)$.

\endroster
\endproclaim
\proclaim{Lemma 1} If \cM\  is a factor, then for any pair of non-zero elements $a, b\in \sM, a, b\neq 0$, then 
$$
a\sM b\neq \lrbrace{0}.
$$
\endproclaim
\demo{Proof} If $a\neq 0$, then 
$$
\sM a \sM=\lrbrace{\sum \xdi a \ydi: \xdone, \cdots, \xdn, \ydone, \cdots \ydn \in \sM}
$$
is a nonzero ideal of $\sM$. Consequently, it is $\sig$-strongly dense in $\sM$.
Hence $\sM a \sM b\neq \lrbrace{0}$, thus $a\sM b\neq \lrbrace{0}$.\QED

\enddemo
\proclaim{Lemma 2} If $\sM$ is a factor, then for every pair of non-zero operators $a, b\in \sM_+$, there exists a non-zero $x\in\sM$ such that
$$
x^*x\leq a, \quad
\text{and }\quad  xx^*\leq b.
$$
\endproclaim
\demo{Proof} As $\sM$ is a factor, the last lemma yields
$$
a^\frac12\sM b^\frac12\neq \lrbrace{0}.
$$
Choose a nonzero $u=b^\frac12 y a^\frac12\in \sM $. Then we get
$$\aligned
u^*u&=a^\frac12 y^*b^\frac12b^\frac12 y a^\frac12\leq \lrnormsq{b^\frac12y}a;\\
uu^*&=b^\frac12 ya^\frac12 a^\frac12 y^* b^\frac12\leq \lrnormsq{ya^\frac12}b.
\endaligned
$$
With 
$$
x=\frac1{\max\lrbrace{\lrnorm{b^\frac12 y}, \lrnorm{ya^\frac12}}}u
$$
we have
$$
x^*x \leq a, \quad xx^* \leq b. \quad x\neq 0.
$$
This completes the proof. \QED
\enddemo

\proclaim{Lemma 3} Fix a $\sig$-finite factor $\sM$ and a nonzero element $x\in \sM$. Set
$$
h=x^*x\quad\text{and}\quad k=xx^*.
$$
Then there exist sequences 
$\lrbrace{\hdn: n\in \N}$, $\lrbrace{\kdn: n\in \N}$ of positive operators in \cM\ and $\lrbrace{\wdn: n\in \N}$ of unitaries in \cM\ such that 
$$
\wdn \hdn\wdnus =\kdn\quad \text{and}\quad h=\sumdnonetoinf \hdn,\quad k=\sumdnonetoinf \kdn.
$$
\endproclaim
\demo{Proof} Let
$$
x=uh^\frac12
$$
be the polar decomposition of $x$. Then we have
$$
k=uhu^*.
$$
Set $e=u^*u$ and $f=uu^*$. Then $e$ and $f$ are  the range projections of $h$ and $k$ respectively. If 
$1-e\sim 1-f$, then choose a partial isometery $v
\in\sM$ such that
$$
v^*v=1-e
\quad\text{and}\quad vv^*=1-f,
$$
and set
$$
w=u+v
$$
to obtain a unitary $w$ such that
$$
k=whw^*.
$$
So in this case, the singleton systems $\lrbrace{h}$, $
\lrbrace{k}$ and $\lrbrace{w}$ answer the question. 

If $\sM$ is finite, then the normalized trace $
\tau$ on $\sM$ regulates the equivalence of projections in 
\cMp So the simple computations:
$$
\gathered
\tau\lr{e}=\tau\lr{u^*u}=
\tau\lr{uu^*}=\tau(f),\\
\tau(1-e)=1-\tau(e)=1-\tau(f)=\tau(1-f),\\
1-e\sim 1-f,
\endgathered
$$
imply our assertion.  
Now, we  assume that $\sM$ is properly infinite. First we set up the notations before proceeding to the next step. Let $\sA$  be a maximal abelian von Neumann subalgebras of \cM containing $h$ and $\sB$ be a maximal abelian subalgebra of \cM which contains $u\sA u^*$. We then consider the semi-finite case and the purely infinite case separately. 

The case that \cM is semi-finite: Let $\tau$ be a faithful semi-finite normal trace on \cM and $\fmdtau$ be the definition ideal of $\tau$, i.e.,
$$
\fmdtau=L^1(\sM, \tau)\cap \sM=\lrbrace{x\in\sM: \tau(\lrabs{x})<+\infty}.
$$
Let $\edz$  be the projection of $\sA$  such that
$$\gathered
\edz\sA=\overline{\fmdtau\cap \sA}^{ \text{weak clouse}}.
\endgathered
$$
Suppose that $\fdz$ is the projection in $\sB$ such that 
$$
\fdz\sB=\overline{\fmdtau\cap \sB}^{ \text{weak clouse}}.
$$
Since the projection $u\edz u^*$ is semi-finite, we have 
$$
u\edz u^*\leq u e\edz u^*\leq f\fdz.
$$
Let $\lrbrace{\edi: i\in I}$ be a maximal orthogonal family of projections in $\sA$ such that $\tau\lr{\edi}<+\infty$. Then $\edz=\sumd{i\in I} \edi$ and hence $e\edz=\sumd{i\in I}e\edi$. Set $\udi=u\edi$. Then we have
$$
\udius\udi=e\edi \quad \text{and}\quad \udi\udius =fu\edi u^*.
$$
The finiteness $\tau\lr{e\edi}<+\infty, i\in I, $ implies that there exists a unitary $\wdi\in\sU(\sM)$ such that
$$
\udi=\wdi e\edi, \quad i\in I.
$$
Hence with $\hdi=h\edi$, we have
$$\gathered
h\edz=\sumd{i\in I}\hdi, \\
\sumd{i\in I}\wdi \hdi \wdius= ku\edz u^*.
\endgathered
$$
Since $e-e\edz$ does not majorize non-zero $\tau$-finite projection at all, $k-ku\edz u^*$ does not subordinate nonzero $\tau$-finite projection either. Since 
$$\gathered
\lr{x\lr{1-\edz}}^*\lr{x\lr{1-\edz}}=h\lr{1-\edz },\\
\lr{x\lr{1-\edz}}\lr{x\lr{1-\edz}}^*=k\lr{1-u\edz u^*}.
\endgathered
$$ 
Restricting our attention to $h\lr{e-e\edz}$ and $k\lr{f-fu\edz u^*}$,	we assume that $\edz=0$. This means that the projections $e$ and $f$ are both purely infinite and the the restrictions of the trace $\tau$ to $\sA e$ and $\sB f$ are both purely infinite.  Since $e\in A$ and the reduced abelian subalgebra $\sA_e$ is maximally abelian in the reduced factor $\sM_e$. Since $\tau(e)=+\infty$, $\sM_e$ is of infinite dimension. Being maximally abelian in $\sM_e$, $\sA e$ contains infinitely many orthogonal projection, i.e., the projection $e$ is the sum of infinitely many  projections $\lrbrace{\edn: n\in \N}\subset \Proj(\sA)$ such that $\tau\lr{\edn}=+\infty$, i.e., 
$$
e=\sumdnonetoinf \edn, \quad \edn\in \Proj(\sA), \quad  	\tau\lr{\edn}=+\infty.
$$
Setting
$$\gathered
\fdn=u\edn u^*, \quad \hdn=h\edn, \\
\kdn=u\hdn u^*,
\endgathered
$$
we obatain 
$$
h=\sumdnonetoinf \hdn, \quad k=\sumdnonetoinf \kdn, \quad 
u\hdn u^*=\kdn, \quad n\in \N.
$$
Since 
$$\gathered
\tau\lr{1-\edn}\geq \tau\lr{\edm}=+\infty, \\
\tau\lr{1-\fdn}\geq \tau\lr{\fdm}=+\infty,
\endgathered\quad n\neq m,
$$
we conclude that
$$
1-\edn\sim 1-\fdn, \quad n\in\N.
$$
Hence for each $n\in\N$, there exists a unitary $\wdn\in\sU(\sM)$ such that
$$
u\edn=\wdn \edn \quad \text{and}\quad \wdn \hdn\wdnus=\kdn.
$$

The case that $\sM$ is purely infinite: In this case, every pair of nonzero projections are equivalent. So the only case we have to consider is the case that $e=1$ and $f\neq 1$ by symmetry. Let $\sA$ and $\sB$ be the maximal abelian subalgebras such that $\sA$ contains $h$ and $\sB$ contains $u\sA u^*$. As in the previous case, $\sA e$ contains a sequence $\lrbrace{\edn: n\in\N}$ such that 
$$
e=\sumdnonetoinf \edn, \quad \edn\neq 0. 
$$
Set 
$$
\fdn=u\edn u^*\in\sB, \quad n\in\N.
$$
Then we have
$$\gathered
1-\edn \geq \edm\neq 0, \\ 1-\fdn \geq \fdm \neq 0,
\endgathered\quad\text{for } m\neq n, 
$$
so that
$$
1-\edn\sim 1-\fdn.
$$
Thus the arguments in the last paragraph shows the existence of a sequence $\lrbrace{\wdn: n\in\N}$ such that
$$\gathered
u\edn=\wdn \edn, \quad \wdn \hdn \wdnus =k\fdn=\kdn,\\
h=\sumdnonetoinf \hdn, \quad k=\sumdnonetoinf \kdn.
\endgathered
$$
This completes the proof. \QED

\enddemo

\demo{\rm Proof of Main Theorem} $\text{(i)}\Rightarrow \text{(ii)}$: Suppose
$\tau(h)\leq \tau(k)$. 
Let $\fX$ be the family of systems $\lrbrace{\xdi: \xdi\in \sM, i\in I}$ such that
$$\gathered
\sumd{i\in I} \xdius \xdi \leq h, \quad \sumd{i\in I} \xdi\xdius \leq k.
\endgathered
$$
Then the set $\fX$ is an inductive set relative to the inclusion ordering. Zorn's lemma yields the existence of a maximal system $\lrbrace{\xdi: i\in I}$ in $\fX$.
If 
$$a=h-\sumd{i\in I}\xdius \xdi \neq 0\quad \text{and}\quad b=k-\sumd{i\in I}\xdi \xdius\neq 0,
$$
then the last lemma implies the existence of a nonzero $\xdz\in \a\sM b$ such that
$$
x_0^*x_0 \leq a=h-\sumd{i\in I}\xdius \xdi \neq 0\quad \text{and}\quad \xdz\xdzus \leq b=k-\sumd{i\in I}\xdi\xdius.
$$
Thus $\lrbrace{\xdz}\cup\lrbrace{\xdi: i\in I}\in \fX$ is a system properly larger than $\lrbrace{\xdi: i\in I}$, violating the maximality of $\lrbrace{\xdi: i\in I}$ in $\fX$. Hence we have either $a=0$ or $b=0$. If $b=0$, then the normality of the trace $\tau$ implies the following:
$$\aligned
\tau(k)&=\tau\lr{\sumd{i\in I}\xdi \xdius}=\sumd{i\in I}\tau\lr{\xdi\xdius}
=\sumd{i\in I}\tau\lr{\xdius\xdi}\\
&=\tau\lr{\sumd{i\in I}\xdius\xdi}
=\tau(h-a)\leq \tau(h).
\endaligned
$$
This computation shows that $a=0$. Thus we have 
$$
h=\sumd{i\in I}\xdius\xdi\quad\text{and}\quad \sumd{i\in I}\xdi\xdius \leq k.
$$
The $\tau$-integrability assumption on $h$:
$$
\infty>\tau(h)=\sumd{i\in I}\tau\lr{\xdius\xdi}
$$ 
shows that the system $\lrbrace{\xdi: i\in I}$ must be countable. The positivity of each term $\tau\lr{\xdius \xdi}$ yields that the enumeration of the index set $I$ does not afect the summation.

(ii)$\Rightarrow$(iii): It is enough to show that if
$$
h=x^*x\quad \text{and}\quad k=xx^*, \quad x\neq 0,
$$
then there exist a sequence $\lrbrace{\hdn: n\in \N}$ of positive operators in $\sM$ and a sequence $\lrbrace{\udn: n\in\N}$ of unitaries in $\sM$ such that 
$$
h=\sumdnonetoinf \hdi\quad \text{and}\quad
k=\sumdnonetoinf \udi\hdi\udius.
$$
In the case that $\tau(1)<+\infty$, every partial isometry $u$ in $\sM$ admits a unitary extension $v\in\sU(\sM)$ so that choosing $u$ to be the phase of $x$, i.e, the partial isometry appearing in the polar decompostion of $x$:
$$
x=uh^\frac12
$$
we get the unitary equivalence:
$$
k=vhv^*.
$$
So we assume that $\tau(1)=+\infty$. Let
$$
h=\int_{\R_+}\la \txd e(\la), \quad x=uh^\frac12
$$
be the spectral decomposition of $h$ and the polar decompostion of $x$. Then we have
$$
uhu^*=k.
$$
Let  $\edz$ and $\fdz$ be the support projections of $h$ and $k$ respectively. If $1-\edz\sim 1-\fdz$, then the partial isometery $u$ is the restriction of a unitary $v$ in $\sM$ so that $k=vhv^*$.
 Set
$$
e_n=e\lr{\left(2^{n-1}, 2^n\right]}, \quad n\in \Z.
$$

(ii)$\Rightarrow$(i): Suppose that $\lrbrace{\xdi: i\in I}$ is a system in $\sM$ such that 
$$
h=\sumd{i\in I} \xdius \xdi \quad \text{and}\quad \sumd{i\in I}\xdi\xdius \leq k.
$$
The complete additivity of the trace $\tau$ shows the following:
$$\aligned
\tau(h)&=\tau\lr{\sumd{i\in I} \xdius \xdi }=\sumd{i\in I}\tau\lr{\xdius\xdi}
=\sumd{i\in I}\tau\lr{\xdi\xdius}=\tau\lr{\sumd{i\in I}\xdi\xdius}\\
&\leq \tau(k).
\endaligned
$$
This completes the proof. \QED

\proclaim{Corollary 1} Under the same assumption as in the main theorem, 
any $h\in \sM_+$ admits a partition $\lrbrace{\xdi: i\in I}$ in \cM\ such that
$$
h=\sumd{i\in I}\xdius\xdi\quad \text{and}\quad \tau(h)=\sumd{i\in I}\xdi\xdius.
$$ 
\endproclaim

\enddemo
\head{\S2. Type \threee case }
\endhead
Let $\sM$ be a factor and $\lrbrace{\tM, \tau, \th}$ be tha associated noncommutative flow of weights on $\sM$, so that 
$$\gathered
\tau\scirc \thds =e^{-s}\tau,\\
\sM=\tM^\th, \quad \sM'\cap \tM=\sC=\text{The center of }\tM,\\
\lrbrace{\sC, \R, \th}=\text{The flow of weights on }\sM.
\endgathered
$$
Let $\fM$ be the algebra of all $\tau$-measurable operators affiliated with $\tM$. With 
$$
E_\la=\chi_{[\la, +\infty)}, \quad\la>0,
$$
a densely defined closed operator $x$ affilated with $\tM$ is $\tau$-measurable if and only if
$$\gathered
\lim_{\la \to +\infty}\tau\lr{E_\la\lr{\lrabs{x}}}=0, \\
\text{equivalently}\\
\tau(E_\la(\lrabs{x})<+\infty \quad \text{for some } \la>0.
\endgathered
$$
The predual $\sMds$ is identified with 
$$\aligned
L^1(\sM)&=\lrbrace{x\in\fM: \thds(x)=e^{-s}x}\\
&=\text{The linear span of } \f, \f\in\frSds(\sM),
\endaligned
$$
where we consider $\f\in\sM_*^+$ as an element of $\fM$, indeed $\lrbrace{\f^{\txti t}: t\in \R, \f\in \frS_*^0}$ generates the core $\tM$ of $\sM$ where $ \frSdsuz$ denotes the set of all faithful normal states of \cMc while $\frSds$ means the set of all normal states. The covariant system $\lrbrace{\sA_\f, \R, \th}$   generated by $\lrbrace{\f^{\txti t}: t\in\R}$ is identified with the system $\lrbrace{\Linfr, \rho}$ where $\rho$ is the translation:
$$\gathered
\lr{\rhods(f)}(x)=f(x+s),\\
\f(x)=\explr{-x}, \\
\tau(f)=\intr f(x)\exp(x)\txd x, \\
E_\la(\f)=\chi_{(-\infty, -\log \la]},
\endgathered\quad f\in\Linfr,\ s, x\in\R, \la\in \R_+.
$$ 
As $\tau$ is semifinite on $\sAdf$, there exists a unique $\tau$ preservaing conditional expectation $\sE_\f$ from $\tM$ to $\sAdf$ such that
$$
\tau=\tau\scirc \sE_\f.
$$

\proclaim{Lemma 1} Suppose that $\sM$ is a $\sig$-finite factor. If $\f$ and $\p$ are two nonzero elements in $\sMdsup$, then there exist an element $x\in \sM$ and a constant $C>0$ such that 
$$
0\neq \lrabs{\p^\frac12 x\f^\frac12}\leq C\f \quad\text{and}\quad 
0\neq\lrabs{\f^\frac12x^*\p^\frac12}\leq C\p.
$$
\endproclaim
\demo{Proof} Let $p=s(\f)$ and $q=s(\p)$ be the support projection of $\f$ and $\p$ respectively. Choose $\f'\in \sMdsup$ and $\p'\in\sMdsup$ such that
$$
1-p=s(\f') \quad \text{and}\quad 1-q=s\lr{\p'},
$$
and set
$$
\Bar \f=\f+\f' \quad\text{and}\quad \Bar\p=\p+\p',
$$
so that $\Bar\f$ and $\Bar \p$ are two faithful normal positive linear functional on $\sM$ such that
$$
p\in \sM_{\Bar\f}, \quad q\in\sM_{\Bar\p}, \quad \f=p\Bar\f=\Bar \f p, \quad\text{and}\quad  \p=q\Bar \p=\Bar\p q.
$$
Choose a nonzero element $y\in q\sM p$, then we have
$$
0\neq \p^\frac12 y \f^\frac12 \in L^1(\sM)=\sMds.
$$
Consider the balanced functional $\Bar\rho$:
$$
\Bar\rho=\pmatrix\Bar\p&0\\0&\Bar\f\endpmatrix\in \sN_{\Bar \rho}=\lr{\sM\ot \M_2(\C)}_{\Bar\rho}.
$$
Observe
$$\aligned
\sig^{\Bar\rho}_t\pmatrix0&y\\0&0
\endpmatrix&=
\pmatrix\Bar\p^{\txti t}&0\\0&\Bar\f^{\txti t}\endpmatrix
\pmatrix0&y\\0&0\endpmatrix
\pmatrix\Bar\p^{-\txti t}&0\\0&\Bar\f^{-\txti t}\endpmatrix\\
&=\pmatrix0&\Bar\p^{\txti t}y\Bar\f^{-\txti t}\\0&0\endpmatrix
=\pmatrix0&\sig^{\Bar\p}_t(y)\RN{\Bar\p}{\Bar\f}{t}\\0&0
\endpmatrix\\
&=\pmatrix 0&\RN{\Bar\p}{\Bar\f}{t}\sig^{\Bar\f}_t(y)\\0&0
\endpmatrix\\
&=\pmatrix0&\sig^{\Bar\p, \Bar \f}_t(y)\\0&0\endpmatrix
\endaligned
$$
For each $f\in\loner$, consider th following integral over $\R$
$$\aligned
\sig^{\brho}_f\pmatrix0&y\\0&0
\endpmatrix&=\intr \sig^{\brho}_t\pmatrix0&y\\0&0
\endpmatrix f(t)\txd t\\
&=\pmatrix 0&\intr \sig^{\Bar\p}_t(y)\RN{\p}{\f}{t}f(t)\txd t\\
0&0\endpmatrix\\
&=\pmatrix0&\sig^{\Bar \p, \Bar\f}_f(y)\\0&0\endpmatrix.
\endaligned
$$
Set
$$\gathered
f_1(x)=\frac1{\sqrt{\pi}}\explr{-x^2},\\
f_\e(x)=\frac 1\e f_1\lr{\frac x\e}, \quad \e>0,
\endgathered\quad x\in \R.
$$
Then we have
$$
\limd{\e\searrow 0}\a_{f_\e}(x)=x\quad \text{in the }\sig\text{-strong* convergence}
$$
for any one parameter automorphism group $\lrbrace{\aldt: t\in\R}$ on a \vna \cMc and the map: $t\in\R\mapsto \aldt\lr{\a_{f_\e}(x)}=\a_{\la_t\lr{f_\e}}(x)$ has an entire extenson to the complex plane $\C$ for every $x\in\sM$.  Consequently for sufficiently small $\e>0$ 
$$\gathered
\sig^{\Bar\rho}_{f_\e}\pmatrix0&y\\0&0\endpmatrix
=\pmatrix0&\sig^{\Bar\p,\Bar\f}_{f_\e}(y)\\0&0\endpmatrix\neq 0.\\
\endgathered
$$
Hence with $x(\e)=\sig^{\Bar\p,\Bar\f}_{f_\e}(y)\neq 0, \e>0$, we have
$$\gathered
qx(\e)p=x(\e),\\
z\in \C\mapsto \sig^{\Bar\p,\Bar\f}_z(x(\e))\in q\sM p
\endgathered
$$
is an entire $\sM$-valued function. Consequently the element
$ \sig^{\Bar\p,\Bar\f}_{-\frac\txti 2}(x(\e))$ is a well-defined element of $q\sM p$. Now we compute in $\fM$ the following:
$$\gathered
\p^\frac12 x(\e)\f^\frac 12=\sig^{\Bar\p, \Bar\f}_{-\frac \txti 2}(x(\e))\f=\p\sig^{\Bar\p, \Bar\f}_{\frac\txti 2}\lr{x(\e)}\\
\aligned
\lrabs{\p^\frac12 x(\e)\f^\frac 12}&=\lr{\lr{\p^\frac12 x(\e)\f^\frac 12}^*\lr{\p^\frac12 x(\e)\f^\frac 12}}^\frac12\\
&=\lr{\lr{\sig^{\Bar\p, \Bar\f}_{-\frac \txti 2}(x(\e))\f}^*\lr{\sig^{\Bar\p, \Bar\f}_{-\frac \txti 2}(x(\e))\f}}^\frac12\\
&\leq \lrnorm{\sig^{\Bar\p, \Bar\f}_{-\frac \txti 2}(x(\e))}\f.
\endaligned
\endgathered
$$
Similarly, we have
$$
\lrabs{\lr{\p^\frac12 x(\e)\f^\frac 12}^*}\leq \lrnorm{\sig^{\Bar\f, \Bar\p}_{-\frac\txti 2}(x(\e)^*)}\p.
$$
Thus with $C=\max\lrbrace{\lrnorm{\sig^{\Bar\p, \Bar\f}_{-\frac \txti 2}(x(\e))}, \lrnorm{\sig^{\Bar\f, \Bar\p}_{-\frac\txti 2}(x(\e)^*)}}$ and $x=x(\e)$ we get the desired inequality. 
\QED

\proclaim{Lemma 2} Consider a pari $\f, \p\in\sMdsup$ of nonzero elements. If $u\in\sM$ is a partial isometry such that
$$
u^*u=s(\f), \quad uu^*=s(\p), \quad \p(x)=\f\lr{uxu^*}, x\in\sM,
$$
Then there exist sequences $\lrbrace{\f_n: n\in\N}, \lrbrace{\p_n: n\in\N}$ of elements in $\sMdsup$ and a sequence of unitaries
$\lrbrace{\udn: n\in\N}\i \sU(\sM)$ such that
$$
\f=\sumdnonetoinf \f_n, \quad \p=\sumdnonetoinf \p_n, \quad \p_n=\f_n\scirc \Ad\lr{\udn}.
$$ 
\endproclaim
\demo{}
\enddemo

\proclaim{Lemma 3} For any pair $\f, \p\in\sMdsup$ of nonzero element, there exists a pair $\f_1, \p_1\in\sMdsup$ of nonzero elements and an inner automorphism $\a\in\Int(\sM)$ such that
$$
\f_1\leq \f, \quad \p_1\leq \p, \quad\text{and}\quad \p_1=\f_1\scirc \a.
$$
\endproclaim
\demo{Proof} Let $x$ and $C>0$ be those appearing in Lemma 1 and set
$$
\f_1=\frac 1C\lrabs{\p^\frac12 x\f^\frac 12}, \quad
\p_1=\frac1C\lrabs{\lr{\p^\frac12 x\f^\frac 12}^*}.
$$
\enddemo

$$\aligned
\intr \explr{-x^2}\txd x&=2\int_{\R_+} \explr{-x^2}\txd x
=2 \sqrt{\int_{\R_+\times \R_+}\explr{-\lr{x^2+y^2}}\txd x\txd y}\\
&=2\sqrt{\int_0^\frac \pi2\int_0^\infty \explr{-r^2}r\txd r\txd \th}\\
&=2 \sqrt{\frac12 \int_0^\frac\pi2\int_0^\infty \explr{-s}\txd s \txd\th}\\
&=2\sqrt{\frac \pi4}=\sqrt {\pi}.
\endaligned
$$

Set
$$
\rho=\p^\frac12 x\f^\frac12\in L^1(\sM)=\sM_*.
$$

let 
$$
\rho=u \omega, \quad \omega=\lrabs{\p^\frac12 x\f^\frac12}=u^*\rho=u^*\lr{\p^\frac12 x\f^\frac12}.
$$
be the polar decompostion of $\rho$.
\enddemo

$$\aligned
\rho&=\sig^\p_{-\frac \txti 2}\lr{x}\p^\frac12\f^\frac12
=\sig^\p_{-\frac \txti 2}\lr{x}\p^\frac12 \f^{-\frac12}\f\\
&=\sig^\p_{-\frac \txti 2}\lr{x}\RN{\p}{\f}{-{\frac\txti 2}}\f.
\endaligned
$$

\enddocument

\input 0-macros.tex
\TagsOnRight 
\topmatter
\title{Noncommutative Integration }
\endtitle

\author{Masamichi Takesaki}
\endauthor

\dedicatory{Dedicated to the Memory of William B. Arveson.}
\enddedicatory
\leftheadtext{Noncommutative Integration}
\rightheadtext{Noncommutative Integration}

\endtopmatter
\document

\proclaim{Main Theorem} Let $\sM$ be a semi-finite factor equipped  with a faithful semi-finite normal trace $\tau$. For two non-negative $\tau$-integrable operators $h$ and $k$ affiliated to $\sM$, the following conditions are equivalent\:
\roster
\item "i)" 
$$
\tau(h)\leq \tau(k),
$$
\item"ii)" There exists a sequence $\lrbrace{\xdi: i\in \N}$ of elements in $\sM$ such that
$$
h=\sum_{i=1}^\infty\xdius \xdi\quad\text{and}\quad \quad \sum_{i=1}^\infty\xdi\xdius \leq k.
$$
\item"iii)" There exists a sequence $\lrbrace{\hdi: i\in\N}$ of positive operators and a sequence $\lrbrace{\udi: i\in\N}$ of unitaries in $\sM$ such that
$$
h=\sum_{i=1}^\infty \hdi \quad\text{and} \quad 
\sum_{i=1}^\infty \udi \hdi\udius\leq k.
$$
\endroster
\endproclaim
\proclaim{Lemma 1} If \cM\  is a factor, then for any pair of non-zero elements $a, b\in \sM, a, b\neq 0$, then 
$$
a\sM b\neq \lrbrace{0}.
$$
\endproclaim
\demo{Proof} If $a\neq 0$, then 
$$
\sM a \sM=\lrbrace{\sum \xdi a \ydi: \xdone, \cdots, \xdn, \ydone, \cdots \ydn \in \sM}
$$
is a nonzero ideal of $\sM$. Consequently, it is $\sig$-strongly dense in $\sM$.
Hence $\sM a \sM b\neq \lrbrace{0}$, thus $a\sM b\neq \lrbrace{0}$.\QED

\enddemo
\proclaim{Lemma 2} If $\sM$ is a factor, then for every pair of non-zero operators $a, b\in \sM_+$, there exists a non-zero $x\in\sM$ such that
$$
x^*x\leq a, \quad
\text{and }\quad  xx^*\leq b.
$$
\endproclaim
\demo{Proof} As $\sM$ is a factor, the last lemma yields
$$
a^\frac12\sM b^\frac12\neq \lrbrace{0}.
$$
Choose a nonzero $u=b^\frac12 y a^\frac12\in \sM $. Then we get
$$\aligned
u^*u&=a^\frac12 y^*b^\frac12b^\frac12 y a^\frac12\leq \lrnormsq{b^\frac12y}a;\\
uu^*&=b^\frac12 ya^\frac12 a^\frac12 y^* b^\frac12\leq \lrnormsq{ya^\frac12}b.
\endaligned
$$
With 
$$
x=\frac1{\max\lrbrace{\lrnorm{b^\frac12 y}, \lrnorm{ya^\frac12}}}u
$$
we have
$$
x^*x \leq a, \quad xx^* \leq b. \quad x\neq 0.
$$
This completes the proof. \QED
\enddemo

\proclaim{Lemma 3} Fix a factor $\sM$ and a nonzero element $x\in \sM$. Set
$$
h=x^*x\quad\text{and}\quad k=xx^*.
$$
Then there exist sequences 
$\lrbrace{\hdn: n\in \N}$, $\lrbrace{\kdn: n\in \N}$ of positive operators in \cM\ and $\lrbrace{\wdn: n\in \N}$ of unitaries in \cM\ such that 
$$
\wdn \hdn\wdnus =\kdn\quad \text{and}\quad h=\sumdnonetoinf \hdn,\quad k=\sumdnonetoinf \kdn.
$$
\endproclaim
\demo{Proof} Let
$$
x=uh^\frac12
$$
be the polar decomposition of $x$. Then we have
$$
k=uhu^*.
$$
Set $e=u^*u$ and $f=uu^*$. Then $e$ and $f$ are  the range projections of $h$ and $k$ respectively. If 
$1-e\sim 1-f$, then choose a partial isometery $v\in\sM$ such that
$$
v^*v=1-e\quad\text{and}\quad vv^*=1-f,
$$
and set
$$
w=u+v
$$
to obtain a unitary $w$ such that
$$
k=whw^*.
$$
So in this case, the singleton systems $\lrbrace{h}$, $\lrbrace{k}$ and $\lrbrace{w}$ answer the question. 

If $\sM$ is finite, then the normalized trace $\tau$ on $\sM$ regulates the equivalence of projections in \cM. So the simple computations:
$$\gathered
\tau\lr{e}=\tau\lr{u^*u}=\tau\lr{uu^*}=\tau(f),\\
\tau(1-e)=1-\tau(e)=1-\tau(f)=\tau(1-f)\\
1-e\sim 1-f,
\endgathered
$$
imply our assertion.  
Now, we  assume that $\sM$ is properly infinite. First we set up the notations before proceeding the next step. Let $\sA$  be a maximal abelian von Neumann subalgebras of \cM containing $h$ and $\sB$ be a maximal abelian subalgebra of \cM which contains $u\sA u^*$. We then consider the semi-finite case and the purely infinite case separately. 

The case that \cM is semi-finite: Let $\tau$ be a faithful semi-finite normal trace on \cM and $\fmdtau$ be the definition ideal of $\tau$, i.e.,
$$
\fmdtau=L^1(\sM, \tau)\cap \sM=\lrbrace{x\in\sM: \tau(\lrabs{x})<+\infty}.
$$
Let $\edz$  be the projection of $\sA$  such that
$$\gathered
\edz\sA=\overline{\fmdtau\cap \sA}^{ \text{weak clouse}}.
\endgathered
$$
Suppose that $\fdz$ is the projection in $\sB$ such that 
$$
\fdz\sB=\overline{\fmdtau\cap \sB}^{ \text{weak clouse}}.
$$
Since the projection $u\edz u^*$ is semi-finite, we have 
$$
u\edz u^*\leq u e\edz u^*\leq f\fdz.
$$
Let $\lrbrace{\edi: i\in I}$ be a maximal orthogonal family of projections in $\sA$ such that $\tau\lr{\edi}<+\infty$. Then $\edz=\sumd{i\in I} \edi$ and hence $e\edz=\sumd{i\in I}e\edi$. Set $\udi=u\edi$. Then we have
$$
\udius\udi=e\edi \quad \text{and}\quad \udi\udius =fu\edi u^*.
$$
The finiteness $\tau\lr{e\edi}<+\infty, i\in I, $ implies that there exists a unitary $\wdi\in\sU(\sM)$ such that
$$
\udi=\wdi e\edi, \quad i\in I.
$$
Hence with $\hdi=h\edi$, we have
$$\gathered
h\edz=\sumd{i\in I}\hdi, \\
\sumd{i\in I}\wdi \hdi \wdius= ku\edz u^*.
\endgathered
$$
Since $e-e\edz$ does not majorize non-zero $\tau$-finite projection at all, $k-ku\edz u^*$ does not subordinate nonzero $\tau$-finite projection either. Since 
$$\gathered
\lr{x\lr{1-\edz}}^*\lr{x\lr{1-\edz}}=h\lr{1-\edz },\\
\lr{x\lr{1-\edz}}\lr{x\lr{1-\edz}}^*=k\lr{1-u\edz u^*}.
\endgathered
$$ 
Restricting our attention to $h\lr{e-e\edz}$ and $k\lr{f-fu\edz u^*}$,	we assume that $\edz=0$. This means that the projections $e$ and $f$ are both purely infinite and the @restriction of the trace $\tau$ on $\sA e$ and $\sB f$ are both purely infinite.

\enddemo

\demo{\rm Proof of Main Theorem} $\text{(i)}\Rightarrow \text{(ii)}$: Suppose
$\tau(h)\leq \tau(k)$. 
Let $\fX$ be the family of systems $\lrbrace{\xdi: \xdi\in \sM, i\in I}$ such that
$$\gathered
\sumd{i\in I} \xdius \xdi \leq h, \quad \sumd{i\in I} \xdi\xdius \leq k.
\endgathered
$$
Then the set $\fX$ is an inductive set relative to the inclusion ordering. Zorn's lemma yields the existence of a maximal system $\lrbrace{\xdi: i\in I}$ in $\fX$.
If 
$$a=h-\sumd{i\in I}\xdius \xdi \neq 0\quad \text{and}\quad b=k-\sumd{i\in I}\xdi \xdius\neq 0,
$$
then the last lemma implies the existence of a nonzero $\xdz\in \sM b$ such that
$$
x_0^*x_0 \leq a=h-\sumd{i\in I}\xdius \xdi \neq 0\quad \text{and}\quad \xdz\xdzus \leq b=k-\sumd{i\in I}\xdi\xdius.
$$
Thus $\lrbrace{\xdz}\cup\lrbrace{\xdi: i\in I}\in \fX$ is a system properly larger than $\lrbrace{\xdi: i\in I}$, violating the maximality of $\lrbrace{\xdi: i\in I}$ in $\fX$. Hence we have either $a=0$ or $b=0$. If $b=0$, then the normality of the trace $\tau$ implies the following:
$$\aligned
\tau(k)&=\tau\lr{\sumd{i\in I}\xdi \xdius}=\sumd{i\in I}\tau\lr{\xdi\xdius}
=\sumd{i\in I}\tau\lr{\xdius\xdi}\\&=\tau\lr{\sumd{i\in I}\xdius\xdi}
=\tau(h-a)\leq \tau(h).
\endaligned
$$
This computation shows that $a=0$. Thus we have 
$$
h=\sumd{i\in I}\xdius\xdi\quad\text{and}\quad \sumd{i\in I}\xdi\xdius \leq k.
$$
The $\tau$-integrability assumption on $h$:
$$
\infty>\tau(h)=\sumd{i\in I}\tau\lr{\xdius\xdi}
$$ 
shows that the system $\lrbrace{\xdi: i\in I}$ must be countable. The positivity of each term $\tau\lr{\xdius \xdi}$ yields that the enumeration of the index set $I$ does not afect the summation.

(ii)$\Rightarrow$(iii): It is enough to show that if
$$
h=x^*x\quad \text{and}\quad k=xx^*, \quad x\neq 0,
$$
then there exist a sequence $\lrbrace{\hdn: n\in \N}$ of positive operators in $\sM$ and a sequence $\lrbrace{\udn: n\in\N}$ of unitaries in $\sM$ such that 
$$
h=\sumdnonetoinf \hdi\quad \text{and}\quad
k=\sumdnonetoinf \udi\hdi\udius.
$$
In the case that $\tau(1)<+\infty$, every partial isometry $u$ in $\sM$ admits a unitary extension $v\in\sU(\sM)$ so that choosing $u$ to be the phase of $x$, i.e, the partial isometry appearing in the polar decompostion of $x$:
$$
x=uh^\frac12
$$
we get the unitary equivalence:
$$
k=vhv^*.
$$
So we assume that $\tau(1)=+\infty$. Let
$$
h=\int_{\R_+}\la \txd e(\la), \quad x=uh^\frac12
$$
be the spectral decomposition of $h$ and the polar decompostion of $x$. Then we have
$$
uhu^*=k.
$$
Let  $\edz$ and $\fdz$ be the support projections of $h$ and $k$ respectively. If $1-\edz\sim 1-\fdz$, then the partial isometery $u$ is the restriction of a unitary $v$ in $\sM$ so that $k=vhv^*$.
 Set
$$
e_n=e\lr{\left(2^{n-1}, 2^n\right]}, \quad n\in \Z.
$$

(ii)$\Rightarrow$(i): Suppose that $\lrbrace{\xdi: i\in I}$ is a system in $\sM$ such that 
$$
h=\sumd{i\in I} \xdius \xdi \quad \text{and}\quad \sumd{i\in I}\xdi\xdius \leq k.
$$
The complete additivity of the trace $\tau$ shows the following:
$$\aligned
\tau(h)&=\tau\lr{\sumd{i\in I} \xdius \xdi }=\sumd{i\in I}\tau\lr{\xdius\xdi}
=\sumd{i\in I}\tau\lr{\xdi\xdius}=\tau\lr{\sumd{i\in I}\xdi\xdius}\\
&\leq \tau(k).
\endaligned
$$
This completes the proof. \QED

\proclaim{Corollary 1} Under the same assumption as in the main theorem, 
any $h\in \sM_+$ admits a partition $\lrbrace{\xdi: i\in I}$ in \cM\ such that
$$
h=\sumd{i\in I}\xdius\xdi\quad \text{and}\quad \tau(h)=\sumd{i\in I}\xdi\xdius.
$$ 
\endproclaim

\enddemo

\enddocument